\documentclass{amsart}
\usepackage{amsmath,amssymb,amsthm}
\usepackage{aliascnt}
\usepackage{cleveref}
\usepackage{setspace}
\usepackage{url}
\allowdisplaybreaks
\newcommand{\newaliasthm}[3]{%
  \newaliascnt{#1}{theorem}%
  \newtheorem{#1}[#1]{#2}%
  \aliascntresetthe{#1}%
  \crefname{#1}{#2}{#3}%
  \Crefname{#1}{#2}{#3}%
}

\theoremstyle{plain}
\newtheorem{theorem}{Theorem}[section]
\crefname{theorem}{Theorem}{Theorems}
\Crefname{theorem}{Theorem}{Theorems}
\newaliasthm{proposition}{Proposition}{Propositions}
\newaliasthm{lemma}{Lemma}{Lemmas}
\newaliasthm{corollary}{Corollary}{Corollaries}
\theoremstyle{definition}
\newaliasthm{definition}{Definition}{Definitions}
\newaliasthm{problem}{Problem}{Problems}
\newaliasthm{example}{Example}{Examples}
\theoremstyle{remark}
\newaliasthm{remark}{Remark}{Remarks}
\newaliasthm{conjecture}{Conjecture}{Conjectures}

\crefname{section}{Section}{Sections}
\Crefname{section}{Section}{Sections}
\crefname{equation}{equation}{equations}
\Crefname{equation}{Equation}{Equations}
\numberwithin{equation}{section}

\title[$(r,s)$-Designs in the Nonbinary Johnson Scheme]{$(r,s)$-Designs in the Nonbinary Johnson Scheme: Constructions and Linear Codes}

\author[R.~Misawa]{Ryutaro Misawa}
\address[R.~Misawa]{Division of Mathematics\\
Graduate School of Information Sciences\\
Tohoku University\\
6-3-09 Aramaki-Aza-Aoba, Aoba-ku, Sendai 980-8579\\
Japan}
\email{misawa.ryutaro.q2@dc.tohoku.ac.jp}

\author[Y.~Tanaka]{Yuuho Tanaka}
\address[Y.~Tanaka]{Mathematical Sciences Program\\
Department of Science and Technology\\
Faculty of Science and Technology\\
Oita University\\
700 Dannoharu, Oita 870-1192\\
Japan}
\email{tanaka-yuuho@oita-u.ac.jp}

\date{}

\begin{document}
\begin{abstract}
Using factorizations of complete graphs, we determine the existence conditions for $(2,1)$-designs in the nonbinary Johnson scheme where each $3$-element subset appears exactly once as a support.
Furthermore, by utilizing large sets of Steiner triple systems, we construct designs with index $1$ that cannot be obtained from the known construction combining block designs and orthogonal arrays.
For linear codes, we express multiplicative character sums over specified nonzero coordinates as sums over the codewords of the dual code.
Using this identity along with the rank of a matrix formed by Krawtchouk polynomials, we provide a sufficient condition for the set of all codewords of a given weight to form an $(r,s)$-design.
In this matrix, we can utilize not only the rows corresponding to weights absent from the code, but also the rows corresponding to weights for which the set of codewords already forms a design.
Moreover, let $q>2$ be a prime power, and let $0 \neq C \subsetneq \mathbb{F}_q^n$ be a linear code with exactly two nonzero weights.
For $t \geq 2$, the set of all codewords of weight $w \geq t$ forms a non-empty $(t,t)$-design if and only if the dual minimum distance is at least $t+2$.
Consequently, this provides a negative answer to the Wu--Zhou problem concerning two-weight codes with a dual minimum distance of $3$.
\end{abstract}
\maketitle

%%%%%%%%%%%%%%%%%%%%%%%%%%%%%%%%%%%%%%%%%%%%%%%%%%%%%%%%%%%%%%%%%%%%%%%
\section{Introduction}\label{sec:introduction}

Let $n \geq 1$, $q \geq 2$, and $0 \leq w \leq n$ be integers, and let $\Omega$ be a $q$-element set with a specified element denoted by $0$.
Set $[n] := \{1, \ldots, n\}$.
For $y = (y_i)_{i \in [n]} \in \Omega^n$, the set of coordinates of its nonzero components is called the support, and its cardinality is called the weight.
Let $Y \subseteq \Omega^n$ be a set of vectors of weight $w$.
For $0 \leq s \leq r \leq w$, $Y$ is called an $(r,s)$-$(n,w,q,\lambda)$ design if for any $S \subseteq R \subseteq [n]$ with $|S|=s$ and $|R|=r$, and any $\omega \in (\Omega \setminus \{0\})^S$,
\[
\left|\{y\in Y\mid y_i\neq0\ (i\in R),\ y_i=\omega_i\ (i\in S)\}\right|=\lambda
\]
holds.
This common non-negative integer $\lambda$ is called the index.
Here, $(\Omega \setminus \{0\})^S$ denotes the set of all mappings from $S$ to $\Omega \setminus \{0\}$.
That is, the number of vectors that are nonzero at $r$ specified coordinates and take specified values at $s$ of those coordinates is independent of the choice of coordinates and values.
For $q = 2$, this reduces to the condition of an $r$-design, where each $r$-element subset is contained in exactly $\lambda$ blocks.
For $w = n$, this becomes the condition of an orthogonal array of alphabet size $q-1$ and strength $s$.

Delsarte \cite{Delsarte1973} provided a unified treatment of block designs and orthogonal arrays through design theory on association schemes.
Designs in the nonbinary Johnson scheme \cite{TAG1985} on vectors of fixed length and weight yield block designs when $q=2$, and orthogonal arrays when $w=n$.
Nozaki and Watanabe \cite[Sections~3--4]{NW2025} combinatorially characterized designs in this scheme and introduced the $(r,s)$-design defined above.
Its number of points is
\[
|Y|=(q-1)^s\frac{\binom nr}{\binom wr}\lambda
\]
\cite[Corollary~4.4]{NW2025}.
Therefore, for fixed parameters $n, w, q, r, s$, a design with index $1$ has the minimum number of points among non-empty designs.

Construction A by Nozaki and Watanabe \cite[Section~6]{NW2025} is a method that places the rows of an orthogonal array on each block of a block design.
Accounting for the multiplicities of blocks and rows, the number of vectors that are nonzero at specified coordinates and take specified values at a subset of those coordinates is equal to the product of the indices of the block design and the orthogonal array used.
Problem~6.3 in \cite{NW2025} asks for the construction of an $(r,s)$-design with index $1$ that cannot be obtained from Construction A.
In this paper, we consider this problem for $(2,1)$-designs of weight $3$.
Furthermore, we investigate the conditions under which codewords of a linear code of a given weight form an $(r,s)$-design.

First, for $n \geq 3$, $q \geq 2$, and $\lambda \geq 1$, we show that the existence conditions for a $(2,1)$-$(n,3,q,\lambda)$ design in which each $3$-element subset appears exactly once as a support are
\[
\lambda(q-1)=n-2,
\qquad \lambda(n-1)\equiv0\pmod2
\]
(Theorem~\ref{thm:all-triples}).
By fixing one coordinate and its nonzero value, and regarding the remaining two nonzero coordinates as an edge, we find that the design condition is equivalent to the resulting graph being $\lambda$-regular.
Thus, we construct the design by factorizing the complete graph into $\lambda$-regular spanning subgraphs.
For $\lambda = 1$ and $n \geq 5$, this coincides with Etzion's construction of constant-weight codes \cite[Theorem~17]{Etzion2022}.

Next, we use a Steiner triple system, a family of $3$-element subsets containing each $2$-element subset exactly once.
A partition of all $3$-element subsets into disjoint Steiner triple systems is called a large set.
From its existence theorem \cite{Ji2005} and the construction by Chee and Ling \cite[Sections~II--III]{CheeLing2007}, we obtain a $(2,1)$-$(n,3,q,1)$ design for
\[
n>7,\qquad n\equiv1,3\pmod6,\qquad 2\leq q\leq n-1
\]
(Theorem~\ref{thm:large-set-construction}).
In this construction, distinct vectors have distinct supports.
On the other hand, in a $(2,1)$-design with index $1$ obtained from Construction A, each support appears $q-1$ times (Lemma~\ref{lem:construction-a}).
Therefore, for $q \geq 3$, this design cannot be obtained from Construction A.

The Assmus--Mattson theorem \cite[Theorem~4.2]{AM1969} states that under certain conditions on the minimum distance of a code and the number of nonzero weights in its dual code, the supports of codewords of a given weight form a block design.
In this paper, we use a form of the theorem \cite[Theorem~1.1]{MT2018} that counts the support of each codeword exactly once, allowing for repeated blocks.
To obtain an $(r,s)$-design, it is necessary to show that, in addition to the number of codewords containing each $r$-element subset in their supports, the number of such codewords taking specified nonzero values at $s$ of those coordinates is also independent of the choice of coordinates and values.

For linear codes, let $q$ be a prime power, $\Omega = \mathbb{F}_q$, and $C \subseteq \mathbb{F}_q^n$ be a nonzero linear code.
Let $\mathbb{F}_q^\times := \mathbb{F}_q \setminus \{0\}$, and denote the set of all codewords of weight $w$ by $C_w$.
Let $d$ be the minimum distance of $C$, and $C^\perp$ be its dual code with respect to the standard inner product.
For $s \geq 1$, due to the invariance of $C_w$ under nonzero scalar multiplication, the condition for an $(r,s)$-design can be characterized by the frequencies of component ratios (Theorem~\ref{thm:projective-characterization}).
Furthermore, we convert multiplicative character sums over specified nonzero coordinates into sums over the dual code via the orthogonality of additive characters.

Using this equality, we provide a sufficient condition based on the rank of a matrix consisting of Krawtchouk polynomials (Theorem~\ref{thm:rank}).
The rows of the matrix can be selected from weights $w$ such that $C_w = \varnothing$ and those for which $C_w$ is already known to form an $(r,s)$-design.
As an example, we consider the Hexacode \cite{MOWW1978}, a linear code over $\mathbb{F}_4$ of length $6$, dimension $3$, and minimum distance $4$.
We show that $C_6$ is a $(3,2)$-design, and then apply the rank condition to deduce that $C_4$ is a $45$-point $(3,2)$-$(6,4,4,1)$ design (Example~\ref{ex:hexacode}).
By polynomial interpolation, one may replace this rank condition with a condition on the cardinality of an integer set determined by the dual weights (Corollary~\ref{cor:ordinary-weight}).
Furthermore, for $2 \leq r < d$, it is sufficient that the number of nonzero weights in the dual code is at most $(d-r)/(r-1)$ (Corollary~\ref{cor:few-dual-weights}).
However, this last condition cannot hold for $r \geq 3$, and the conclusion for $r=2$ also follows from Wu and Zhou \cite[Theorem~3.6]{WZ2026} (Remark~\ref{rem:new-few-dual-weights}).

Finally, we consider linear codes with exactly two nonzero weights, namely two-weight codes.
A $(t,t)$-design is called a $q$-ary $t$-design.
We denote the weight of a vector $u$ by $\operatorname{wt}(u)$, and for $C^\perp \neq \{0\}$, we write its minimum distance as
\[
d(C^\perp):=\min\{\operatorname{wt}(u)\mid u\in C^\perp\setminus\{0\}\}.
\]
Let $q > 2$, $0 \neq C \subsetneq \mathbb{F}_q^n$, and $t \geq 2$, and suppose that $C$ is a two-weight code.
Then we have
\[
\begin{gathered}
\exists\,w\in\{1,\ldots,n\},\ w\geq t:\quad
C_w\neq\varnothing,\quad C_w\text{ is a }q\text{-ary }t\text{-design}
\\
\Longleftrightarrow\quad d(C^\perp)\geq t+2
\end{gathered}
\]
(Theorem~\ref{thm:new-two-weight-exact}).
The sufficiency is due to Wu and Zhou \cite[Theorem~3.6]{WZ2026}.
To prove necessity, we show that since one $C_w$ is a design, the set of codewords for the other nonzero weight also forms a design.
Next, we fix the values of specified coordinates and apply a MacWilliams-type transform to a polynomial whose exponents are the weights of the remaining components, showing that the set of codewords for each nonzero weight of the dual code also forms a design.
The case where the dual minimum distance is $t+1$ contradicts the minimum distance condition for $t \geq 3$, and the two-weight code condition for $t=2$.
In particular, for two-weight codes with dual minimum distance $3$, a non-empty $C_w$ never forms a $q$-ary $t$-design for $t \geq 2$ (Corollary~\ref{cor:new-WZ-open-problem}).
This gives a negative answer to the open problem raised by Wu and Zhou \cite[p.~16, after Remark~4.15]{WZ2026}.

In Section~\ref{sec:definitions}, we define $(r,s)$-designs and describe Construction A.
In Section~\ref{sec:bounds}, we derive an upper bound on the number of points from the minimum distance of constant-weight codes.
In Section~\ref{sec:weight-three}, we construct designs of weight $3$.
In Section~\ref{sec:linear}, for the set of codewords of a linear code with a given weight, we prove a characterization by component ratios and character sums, sufficient conditions using the dual code, and a necessary and sufficient condition for two-weight codes.

%%%%%%%%%%%%%%%%%%%%%%%%%%%%%%%%%%%%%%%%%%%%%%%%%%%%%%%%%%%%%%%%%%%%%%%
\section{Nonbinary Johnson schemes and $(r,s)$-designs}\label{sec:definitions}

For a finite set $X$ and a non-negative integer $k$, we write $\binom Xk:=\{A\subseteq X\mid |A|=k\}$.
We continue to use $n, q$, and $\Omega$ introduced in Section~\ref{sec:introduction}, and set $\Omega^\times:=\Omega\setminus\{0\}$.
For $A\subseteq[n]$, let $\Omega^A$ denote the set of all mappings from $A$ to $\Omega$.
The direct product indexed by the empty set is defined as $(\Omega^\times)^\varnothing:=\{()\}$, and the empty product is defined to be $1$.
The support $\operatorname{supp}(x)$ and the Hamming weight $\operatorname{wt}(x)$ of $x=(x_i)_{i\in[n]}\in\Omega^n$ are defined as follows:
\[
\operatorname{supp}(x):=\{i\in[n]\mid x_i\neq0\},
\qquad
\operatorname{wt}(x):=|\operatorname{supp}(x)|.
\]
For $A\subseteq[n]$, we write $x|_A:=(x_i)_{i\in A}$, and for $0\leq w\leq n$, we set
\[
W_w(n,q):=\{x\in\Omega^n\mid\operatorname{wt}(x)=w\}.
\]
The Hamming distance $d_H(x,y)$ is defined as $d_H(x,y):=|\{i\in[n]\mid x_i\neq y_i\}|$.

For $q\geq3$, classifying $(x,y)\in W_w(n,q)\times W_w(n,q)$ by the two values
\[
|\{i\in[n]\mid x_i=y_i\neq0\}|
\quad\text{and}\quad
|\operatorname{supp}(x)\cap\operatorname{supp}(y)|
\]
yields a partition of the relations of the nonbinary Johnson scheme \cite{NW2025, TAG1985}.
In the following, we define designs by counting the number of vectors whose entries at specified coordinates are nonzero and some of which take specified values.

For $Y\subseteq W_w(n,q)$, $S\subseteq R\subseteq[n]$, and $\omega\in(\Omega^\times)^S$, we set
\[
m_{R,S}(Y,\omega)
:=|\{y\in Y\mid R\subseteq\operatorname{supp}(y),\ y|_S=\omega\}|.
\]
When $S=\varnothing$, we write $m_R(Y):=|\{y\in Y\mid R\subseteq\operatorname{supp}(y)\}|$.

\begin{definition}\label{def:rs-design}
Let $0\leq s\leq r\leq w\leq n$ be integers, and let $\lambda$ be a non-negative integer.
A subset $Y\subseteq W_w(n,q)$ is called an $(r,s)$-$(n,w,q,\lambda)$ design if 
\[
m_{R,S}(Y,\omega)=\lambda
\]
for any $S\subseteq R\subseteq[n]$ with $|S|=s$ and $|R|=r$, and for any $\omega\in(\Omega^\times)^S$ \cite[Definition~4.1]{NW2025}.
The value $\lambda$ is called the index, and we simply write $(r,s)$-design when the other parameters are omitted.
\end{definition}

We regard the support of each vector as a block.
Since distinct vectors may have the same support, the family of supports is treated as a multiset.
Let $0\leq t\leq k\leq v$ and $\lambda\geq0$ be integers. Let $V$ be a $v$-element set, and let $\mathcal B$ be a finite multiset consisting of elements of $\binom Vk$.
The elements of $\mathcal B$ are called blocks.
If 
\[
|\{B\in\mathcal B\mid T\subseteq B\}|=\lambda
\]
holds for any $T\in\binom Vt$, then $(V,\mathcal B)$ is called a $t$-$(v,k,\lambda)$ design \cite{BJL1999}.
Here, the number of elements on the left-hand side is counted with multiplicities.
When $V$ is clear from the context, $\mathcal B$ itself is also referred to as a $t$-design.

Taking the sum of $m_{R,S}(Y,\omega)$ over all possible nonzero value assignments $\omega\in(\Omega^\times)^S$ yields the number of vectors whose support contains $R$.
Furthermore, by double counting the number of pairs $(y,R)$ where $y\in Y$ and $R$ is an $r$-element subset of $\operatorname{supp}(y)$, we can determine $|Y|$ \cite[Corollary~4.4]{NW2025}.

\begin{lemma}\label{lem:size}
If $Y$ is an $(r,s)$-$(n,w,q,\lambda)$ design, then for any $R\in\binom{[n]}r$,
\[
m_R(Y)=(q-1)^s\lambda,
\qquad
|Y|=(q-1)^s\frac{\binom nr}{\binom wr}\lambda.
\]
\end{lemma}

\begin{proof}
Fix $S\in\binom Rs$.
Summing $m_{R,S}(Y,\omega)=\lambda$ over $\omega\in(\Omega^\times)^S$ yields the first equation.
Counting the elements of the set $\{(y,R)\mid y\in Y,\ R\in\binom{\operatorname{supp}(y)}r\}$ in two ways, we obtain
\[
|Y|\binom wr=\sum_{R\in\binom{[n]}r}m_R(Y)=\binom nr(q-1)^s\lambda,
\]
which implies the equation for the size.
\end{proof}

Let $v\geq1$, $N\geq0$, and $0\leq t\leq k$ be integers.
An $N\times k$ array over a $v$-element set is called an orthogonal array $\operatorname{OA}(N,k,v,t)$ if, for any choice of $t$ columns, every possible tuple of entries in those $t$ columns appears the same number of times across all rows \cite{HSS1999}.
If each such tuple appears $\lambda$ times, the index is $\lambda$, and $N=\lambda v^t$.

Construction A by Nozaki and Watanabe \cite[Section~6]{NW2025} places the rows of an orthogonal array onto each block of a block design.
The case with index $1$ is described below.

\begin{lemma}\label{lem:construction-a}
Let $1\leq s\leq r\leq w\leq n$.
Let $\mathcal B$ be an $r$-$(n,w,1)$ design, and let $A$ be an $\operatorname{OA}((q-1)^s,w,q-1,s)$ with index $1$ over $\Omega^\times$.
For each $B=\{i_1<\cdots<i_w\}\in\mathcal B$ and each row $a=(a_1,\ldots,a_w)$ of $A$, we define $y^{B,a}\in W_w(n,q)$ by
\[
y^{B,a}_{i_j}:=a_j\quad(1\leq j\leq w),
\qquad y^{B,a}_i:=0\quad(i\notin B).
\]
Then
\[
Y:=\{y^{B,a}\mid B\in\mathcal B,\ a\text{ is a row of }A\}
\]
is an $(r,s)$-$(n,w,q,1)$ design, and for each $B\in\mathcal B$,
\[
|\{y\in Y\mid\operatorname{supp}(y)=B\}|=(q-1)^s.
\]
\end{lemma}

\begin{proof}
Fix $R\in\binom{[n]}r$, $S\in\binom Rs$, and $\omega\in(\Omega^\times)^S$.
There exists exactly one $B\in\mathcal B$ containing $R$.
Since the index of the orthogonal array is also $1$, there is exactly one row in $A$ whose entries corresponding to $S$ on $B$ match $\omega$.
Therefore, $m_{R,S}(Y,\omega)=1$.

Since the index is $1$, there are no duplicates among the blocks in $\mathcal B$ or the rows in $A$.
Vectors obtained from distinct blocks have different supports, and distinct rows on the same block also yield distinct vectors.
Thus, the number of vectors having a given $B$ as their support is equal to the number of rows in $A$, which is $(q-1)^s$.
\end{proof}

In Construction A, accounting for the multiplicities of blocks and rows, the number of vectors that are nonzero on $R$ and match $\omega$ on $S$ is equal to the product of the two indices.
If this product is $1$, there are no duplicate blocks or rows, and we obtain the set described in Lemma~\ref{lem:construction-a}.
Nozaki and Watanabe \cite{NW2025} posed the following problem.

\begin{problem}\cite[Problem~6.3]{NW2025}
Construct $(r,s)$-designs with index $1$ which are not obtained from Construction A.
\end{problem}

%%%%%%%%%%%%%%%%%%%%%%%%%%%%%%%%%%%%%%%%%%%%%%%%%%%%%%%%%%%%%%%%%%%%%%%
\section{Constant-weight codes and point bounds}\label{sec:bounds}

If at most one codeword satisfies a given specification of nonzero coordinates and their values on a subset, an upper bound on the code size can be derived.
This section establishes this bound and further demonstrates that if the minimum distance of the code is at least $2w-r-s+1$, no two codewords can simultaneously satisfy the condition.

Let $\mathcal C\subseteq W_w(n,q)$ be a $q$-ary constant-weight code of length $n$ and weight $w$ \cite{Etzion2022}.
For $|\mathcal C|\geq2$, the minimum distance of $\mathcal C$ is defined as follows:
\[
d_H(\mathcal C):=\min\{d_H(x,y)\mid x,y\in\mathcal C,\ x\neq y\}.
\]

\begin{proposition}\label{prop:packing}
Let $0\leq s\leq r\leq w\leq n$ and $\mathcal C\subseteq W_w(n,q)$.
If $m_{R,S}(\mathcal C,\omega)\leq1$ for any $S\subseteq R\subseteq[n]$ with $|S|=s$ and $|R|=r$, and any $\omega\in(\Omega^\times)^S$, then
\[
|\mathcal C|\leq(q-1)^s\frac{\binom nr}{\binom wr}.
\]
Equality holds if and only if $\mathcal C$ is an $(r,s)$-$(n,w,q,1)$ design.
\end{proposition}

\begin{proof}
Counting the number of tuples $(x,R,S)$ satisfying $x\in\mathcal C$, $S\subseteq R\subseteq\operatorname{supp}(x)$, $|R|=r$, and $|S|=s$ in two ways yields
\[
\begin{aligned}
|\mathcal C|\binom wr\binom rs
&=\sum_{R\in\binom{[n]}r}\sum_{S\in\binom Rs}
  \sum_{\omega\in(\Omega^\times)^S}m_{R,S}(\mathcal C,\omega)\\
&\leq\binom nr\binom rs(q-1)^s.
\end{aligned}
\]
By assumption, $m_{R,S}(\mathcal C,\omega)\leq1$ for any $R,S,\omega$, so equality holds in the above inequality if and only if $m_{R,S}(\mathcal C,\omega)=1$ for all such $R,S,\omega$.
This coincides with the condition for an $(r,s)$-design with index $1$ in Definition~\ref{def:rs-design}.
\end{proof}

We evaluate the Hamming distance between two vectors whose supports share $r$ coordinates, and where the values of the nonzero components agree on $s$ of these coordinates.

\begin{proposition}\label{prop:fixed-coordinates}
Let $0\leq s\leq r\leq w\leq n$.
Fix $S\subseteq R\subseteq[n]$ with $|S|=s$ and $|R|=r$, and $\omega\in(\Omega^\times)^S$.
Let
\[
A:=\{x\in W_w(n,q)\mid R\subseteq\operatorname{supp}(x),\ x|_S=\omega\}.
\]
For any $x,y\in A$, we have $d_H(x,y)\leq2w-r-s$ and
\[
|A|=\binom{n-r}{w-r}(q-1)^{w-s}.
\]
Furthermore, for any $\mathcal C$ satisfying the assumption of Proposition~\ref{prop:packing}, the inequality $|\mathcal C|\,|A|\leq|W_w(n,q)|$ holds.
\end{proposition}

\begin{proof}
For $x,y\in W_w(n,q)$, we have
\[
d_H(x,y)=2w-|\operatorname{supp}(x)\cap\operatorname{supp}(y)|
               -|\{i\in[n]\mid x_i=y_i\neq0\}|.
\]
If $x,y\in A$, their supports overlap in at least $r$ coordinates, and their nonzero values agree on at least $s$ coordinates.
Therefore, $d_H(x,y)\leq2w-r-s$.
To construct the elements of the set $A$, we can choose $w-r$ nonzero coordinates from the complement of $R$ and assign nonzero values to the remaining $w-s$ coordinates.
Thus, $|A|=\binom{n-r}{w-r}(q-1)^{w-s}$.
The final inequality is obtained by multiplying both sides of the upper bound on $|\mathcal{C}|$ in Proposition~\ref{prop:packing} by $|A|$ and using the relations
\[
\binom nr\binom{n-r}{w-r}=\binom nw\binom wr,
\text{ and }
|W_w(n,q)|=\binom nw(q-1)^w.
\]
\end{proof}

If the minimum distance of the code strictly exceeds the upper bound on the Hamming distance in Proposition~\ref{prop:fixed-coordinates}, the assumption of Proposition~\ref{prop:packing} is satisfied, allowing the proposition to be applied.
Note that Etzion \cite[Section~4, Lemma~8]{Etzion2022} extended the point bound obtained by Delsarte \cite{Delsarte1973} to nonbinary constant-weight codes.
Substituting the $|A|$ obtained in Proposition~\ref{prop:fixed-coordinates} into the inequality $|\mathcal C|\,|A|\leq |W_w(n,q)|$ yields the same point bound as Proposition~\ref{prop:packing}.

\begin{corollary}\label{cor:distance-bound}
Let $0\leq s\leq r\leq w\leq n$, $\mathcal C\subseteq W_w(n,q)$ with $|\mathcal C|\geq2$, and suppose $d_H(\mathcal C)\geq2w-r-s+1$.
Then
\[
|\mathcal C|\leq(q-1)^s\frac{\binom nr}{\binom wr}.
\]
If equality holds, then $\mathcal C$ is an $(r,s)$-$(n,w,q,1)$ design.
\end{corollary}

\begin{proof}
If there exist distinct $x,y\in\mathcal C$ satisfying the conditions for the same $R,S$ and $\omega$, then Proposition~\ref{prop:fixed-coordinates} implies $d_H(x,y)\leq2w-r-s$, which contradicts the assumption on the minimum distance of $\mathcal C$.
Therefore, $m_{R,S}(\mathcal C,\omega)\leq1$, and the conclusion follows from Proposition~\ref{prop:packing}.
\end{proof}

%%%%%%%%%%%%%%%%%%%%%%%%%%%%%%%%%%%%%%%%%%%%%%%%%%%%%%%%%%%%%%%%%%%%%%%
\section{Constructions of weight $3$}\label{sec:weight-three}

Fixing nonzero values at specified coordinates and removing them from the supports allows the condition for an $(r,s)$-design to be expressed as that of a block design.
It has also been shown by Nozaki--Watanabe \cite[Proposition~4.2]{NW2025} that the sets obtained by removing the coordinates of the specified nonzero components form a design.
To distinguish cases where the same support is obtained from different vectors, we will use multisets in what follows.

\begin{proposition}\label{prop:derived}
Let $0\leq s\leq r\leq w\leq n$, $\lambda\geq0$, and $Y\subseteq W_w(n,q)$.
Then $Y$ is an $(r,s)$-$(n,w,q,\lambda)$ design if and only if for any $S\in\binom{[n]}s$ and $\omega\in(\Omega^\times)^S$, the multiset on $[n]\setminus S$ given by
\[
\bigl(\operatorname{supp}(y)\setminus S\bigr)_{
  y\in Y,\ S\subseteq\operatorname{supp}(y),\ y|_S=\omega}
\]
is an $(r-s)$-$(n-s,w-s,\lambda)$ design.
\end{proposition}

\begin{proof}
Fix $S\in\binom{[n]}s$ and $\omega\in(\Omega^\times)^S$.
The number of blocks containing $T\in\binom{[n]\setminus S}{r-s}$ is equal to $m_{S\cup T,S}(Y,\omega)$.
Since the mapping $T\mapsto S\cup T$ is a bijection from the set of all such $T$ to the set of all $r$-element subsets containing $S$, the equivalence follows.
\end{proof}

For weight $3$, a support excluding a single coordinate can be regarded as an edge of a graph.
In what follows, all graphs are assumed to be finite and simple.
A graph where every vertex has degree $\lambda$ is called $\lambda$-regular, and a $\lambda$-regular subgraph containing all vertices of a given graph is called a $\lambda$-factor.
A partition of the edge set into the edge sets of $\lambda$-factors is called a $\lambda$-factorization.
The conditions for the existence of a $\lambda$-factorization of a complete graph can be derived from Baranyai's factorization theorem \cite{Baranyai1975,Bahmanian2014}.

\begin{lemma}\label{lem:factorization}
Let $N\geq2$ and $\lambda\geq1$ be integers.
A complete graph $K_N$ on $N$ vertices has a $\lambda$-factorization if and only if
\[
\lambda\mid N-1
\text{ and }
\lambda N\equiv0\pmod2.
\]
\end{lemma}

\begin{proof}
Let $k$ be the number of factors.
From the degree of each vertex, we obtain $k\lambda=N-1$.
Furthermore, the two conditions are necessary since the number of edges in each factor is $\lambda N/2$.

Conversely, assume the two conditions hold and let $k=(N-1)/\lambda$.
If $N=2$, then $\lambda=1$, and $K_2$ itself is a $1$-factor.
For $N\geq3$, apply Baranyai's factorization theorem \cite{Baranyai1975,Bahmanian2014} by setting the edge size to $2$, the edge multiplicity to $1$, and the degree of each factor to $\lambda$.
The conditions for this theorem are $2\mid\lambda N$ and $k\lambda=N-1$, both of which are satisfied by our assumptions. 
Therefore, the required factorization exists.
\end{proof}

\begin{proposition}\label{prop:weight-three}
Let $n\geq3$, $q\geq2$, and $\lambda\geq1$.
Suppose that for $Y\subseteq W_3(n,q)$, the mapping $y\mapsto\operatorname{supp}(y)$ is injective.
Then $Y$ is a $(2,1)$-$(n,3,q,\lambda)$ design if and only if for any $a\in[n]$ and $\omega\in\Omega^\times$, the graph with vertex set $[n]\setminus\{a\}$ and edge set
\[
\{\operatorname{supp}(y)\setminus\{a\}\mid
  y\in Y,\ a\in\operatorname{supp}(y),\ y_a=\omega\}
\]
is $\lambda$-regular.
In this case, $\{\operatorname{supp}(y)\mid y\in Y\}$ is a $2$-$(n,3,\lambda(q-1))$ design.
\end{proposition}

\begin{proof}
Substitute $(r,s,w)=(2,1,3)$ into Proposition~\ref{prop:derived}.
Since the supports are distinct, the resulting edges have no multiplicity.
The degree of each $b\in[n]\setminus\{a\}$ is equal to $m_{\{a,b\},\{a\}}(Y,\omega)$, yielding the equivalence in the first part.
Furthermore, the number of supports containing $\{a,b\}$ is the sum of these degrees over all $\omega\in\Omega^\times$, which is $\lambda(q-1)$.
\end{proof}

For odd $n\geq5$, Etzion used a $1$-factorization of $K_{n-1}$ to construct a code of length $n$, weight $3$, alphabet size $n-1$, and minimum distance $4$, in which each $3$-element subset appears exactly once as a support \cite[Theorem~17]{Etzion2022}.
In the following theorem, we determine the existence conditions for any positive integer $\lambda$.

\begin{theorem}\label{thm:all-triples}
Let $n\geq3$, $q\geq2$, and $\lambda\geq1$.
A $(2,1)$-$(n,3,q,\lambda)$ design $Y\subseteq W_3(n,q)$ in which each $3$-element subset appears exactly once as a support exists if and only if
\[
\lambda(q-1)=n-2
\text{ and }
\lambda(n-1)\equiv0\pmod2.
\]
\end{theorem}

\begin{proof}
Suppose there exists a $(2,1)$-$(n,3,q,\lambda)$ design $Y\subseteq W_3(n,q)$ in which each $3$-element subset appears exactly once as a support.
Since each $2$-element subset is contained in exactly $n-2$ supports, Proposition~\ref{prop:weight-three} implies $\lambda(q-1)=n-2$.
Additionally, since the graph in Proposition~\ref{prop:weight-three} is a $\lambda$-regular graph on $n-1$ vertices, the number of edges $\lambda(n-1)/2$ must be an integer.

Conversely, assume the two conditions hold.
By Lemma~\ref{lem:factorization}, for each $a\in[n]$, the complete graph on $[n]\setminus\{a\}$ can be factored into $q-1$ $\lambda$-factors.
Index the factors of each decomposition by $\Omega^\times$ and denote them by $F_a^\omega$.
Let $E(F_a^\omega)$ be its edge set.
For each $B\in\binom{[n]}3$ and $a\in B$, the index $\omega$ of the factor containing $B\setminus\{a\}$ as an edge is uniquely determined.
Using this $\omega$, define
\[
y_a^B:=
\begin{cases}
\omega & a\in B,\ B\setminus\{a\}\in E(F_a^\omega),\\
0 & a\notin B
\end{cases}
\qquad(a\in[n])
\]
and let
\[
Y:=\{y^B\mid B\in\binom{[n]}3\}.
\]
Then $\operatorname{supp}(y^B)=B$, and the graph in Proposition~\ref{prop:weight-three} coincides with $F_a^\omega$.
Therefore, $Y$ is the desired design.
\end{proof}

\begin{corollary}\label{cor:all-triples-index-one}
Let $n\geq3$ and $q\geq2$.
A $(2,1)$-$(n,3,q,1)$ design in which each $3$-element subset appears exactly once as a support exists if and only if $q=n-1$ and $n$ is odd.
For $n\geq5$, the construction in Theorem~\ref{thm:all-triples} coincides with Etzion's construction using $1$-factorizations \cite[Theorem~17]{Etzion2022}.
\end{corollary}

Let $n\geq3$.
A $2$-$(n,3,1)$ design is called a Steiner triple system of order $n$.
A partition of $\binom{[n]}3$ into $n-2$ Steiner triple systems is called a large set of order $n$ \cite{CheeLing2007,Ji2005}.
The following existence conditions are known \cite{Ji2005}.

\begin{lemma}\label{lem:large-set}
For $n>7$, a large set of Steiner triple systems of order $n$ exists if and only if $n\equiv1,3\pmod6$.
\end{lemma}

In the construction by Chee and Ling \cite[Sections~II--III]{CheeLing2007}, different nonzero values are assigned to multiple block-disjoint families.
Applying this to a large set yields designs with index $1$ that cannot be obtained from Construction A when $q\geq3$.

\begin{theorem}\label{thm:large-set-construction}
Let $n>7$, $n\equiv1,3\pmod6$, and $2\leq q\leq n-1$.
Then there exists a $(2,1)$-$(n,3,q,1)$ design whose distinct elements have distinct supports.
If $q\geq3$, there exist designs that cannot be obtained from Construction A.
\end{theorem}

\begin{proof}
By Lemma~\ref{lem:large-set}, take a large set
\[
\binom{[n]}3=\mathcal S_1\sqcup\cdots\sqcup\mathcal S_{n-2}.
\]
Since $q-1\leq n-2$, we use $\mathcal S_1,\ldots,\mathcal S_{q-1}$ among these.
Let $\Omega^\times=\{\omega_1,\ldots,\omega_{q-1}\}$, and for each $B\in\mathcal S_j$, $1\leq j\leq q-1$, define
\[
y_i^B:=
\begin{cases}
\omega_j&i\in B,\\
0&i\notin B
\end{cases}
\qquad(i\in[n])
\]
and let
\[
Y:=\{y^B\mid B\in\mathcal S_j,\ 1\leq j\leq q-1\}.
\]
For any distinct $a,b\in[n]$ and any $\omega_j\in\Omega^\times$, there exists exactly one block in $\mathcal S_j$ containing $\{a,b\}$.
Therefore, $m_{\{a,b\},\{a\}}(Y,\omega_j)=1$, and $Y$ is a $(2,1)$-$(n,3,q,1)$ design.

Since the chosen Steiner triple systems are disjoint, the distinct elements of $Y$ have distinct supports.
On the other hand, when $q\geq3$, Lemma~\ref{lem:construction-a} implies that each support used in a $(2,1)$-design with index $1$ obtained from Construction A appears $q-1>1$ times.
Therefore, this $Y$ cannot be obtained from Construction A.
\end{proof}

\begin{example}\label{ex:sts-nine}
Let $\Omega=\{0,1,2\}$, and write $ijk:=\{i,j,k\}$.
In
\[
\begin{aligned}
\mathcal S_1&:=\{123,147,159,168,249,258,267,348,357,369,456,789\},\\
\mathcal S_2&:=\{124,135,167,189,239,256,278,347,368,458,469,579\},
\end{aligned}
\]
any $2$-element subset is contained in exactly one block of each $\mathcal S_j$, and $\mathcal S_1\cap\mathcal S_2=\varnothing$.
Thus, these are disjoint Steiner triple systems of order $9$.
As in the proof of Theorem~\ref{thm:large-set-construction}, we place the value $1$ on $B\in\mathcal S_1$, the value $2$ on $B\in\mathcal S_2$, and $0$ on the other components.
These $24$ vectors form a $(2,1)$-$(9,3,3,1)$ design.
Since all supports are distinct, Lemma~\ref{lem:construction-a} implies that it cannot be obtained from Construction A.
\end{example}

%%%%%%%%%%%%%%%%%%%%%%%%%%%%%%%%%%%%%%%%%%%%%%%%%%%%%%%%%%%%%%%%%%%%%%%
\section{$(r,s)$-Designs from Linear Codes}\label{sec:linear}

The support of a codeword is invariant under scalar multiplication by a nonzero element.
In this section, we express the condition for an $(r,s)$-design using coordinate ratios, and further rewrite the multiplicative character sum into a sum over the codewords of the dual code.

Let $q$ be a prime power and $\Omega=\mathbb F_q$.
A linear subspace of $\mathbb F_q^n$ is called a linear code of length $n$, and its elements are called codewords \cite{MS1977}.
Let $C\subseteq\mathbb F_q^n$ be a nonzero linear code. For $u,v\in\mathbb F_q^n$, we define $u\cdot v:=\sum_{i=1}^n u_iv_i$.
The dual code $C^\perp$ and the minimum distance $d$ are defined as
\[
C^\perp:=\{u\in\mathbb F_q^n\mid u\cdot c=0\ (c\in C)\},
\qquad
d:=\min\{\operatorname{wt}(c)\mid c\in C\setminus\{0\}\}.
\]
A matrix whose rows form a basis of $C$ is called a generator matrix of $C$.
For $1\leq w\leq n$, the set of all codewords of weight $w$ is denoted by
\[
C_w:=\{c\in C\mid\operatorname{wt}(c)=w\}.
\]

%%%%%%%%%%%%%%%%%%%%%%%%%%%%%%%%%%
\subsection{Characterization by coordinate ratios}

Nonzero codewords belonging to the same one-dimensional subspace have the same support and coordinate ratios.
Therefore, instead of specifying the values of nonzero coordinates, one can count one-dimensional subspaces by specifying their coordinate ratios.
For $c\in C_w$, let
\[
\langle c\rangle:=\{ac\mid a\in\mathbb F_q\},
\qquad
\mathcal P_w:=\{\langle c\rangle\mid c\in C_w\}.
\]
Each $\langle c\rangle\in\mathcal P_w$ contains exactly $q-1$ elements of $C_w$.
For $\operatorname{supp}(c)$ and $\varnothing\neq S\subseteq\operatorname{supp}(c)$, the coordinate ratios $c_i/c_{\min S}$ ($i\in S$) do not depend on the choice of the nonzero representative element of $\langle c\rangle$.

\begin{theorem}\label{thm:projective-characterization}
Let $1\leq s\leq r\leq w\leq n$ and an integer $\lambda\geq0$.
Then $C_w$ is an $(r,s)$-$(n,w,q,\lambda)$ design if and only if for any $S\subseteq R\subseteq[n]$ with $|S|=s$ and $|R|=r$, and any $\eta\in(\mathbb F_q^\times)^{S\setminus\{\min S\}}$, we have
\[
\left|\left\{\langle c\rangle\in\mathcal P_w\ \middle|\
R\subseteq\operatorname{supp}(c),\quad
\frac{c_i}{c_{\min S}}=\eta_i\quad(i\in S\setminus\{\min S\})
\right\}\right|=\lambda.
\]
\end{theorem}

\begin{proof}
Fix $S\subseteq R$ and $\omega\in(\mathbb F_q^\times)^S$, and set $\eta_i:=\omega_i/\omega_{\min S}$.
For each $\langle c\rangle$ in the set on the left-hand side of the above equation, there exists exactly one codeword whose coordinates on $S$ coincide with $\omega$.
This codeword is given by $(\omega_{\min S}/c_{\min S})c$, where the scalar is uniquely determined by the coordinate at $\min S$.

Therefore, the cardinality of the above set is equal to $m_{R,S}(C_w,\omega)$.
Conversely, any $\eta$ can be obtained by setting $\omega_{\min S}=1$ and $\omega_i=\eta_i$, which implies that the two conditions are equivalent.
\end{proof}

\begin{corollary}\label{cor:projective-supports}
Let $1\leq s\leq r\leq w\leq n$ and suppose that $C_w$ is an $(r,s)$-$(n,w,q,\lambda)$ design.
Then the multiset of subsets of $[n]$ defined by
\[
\bigl(\operatorname{supp}(c)\bigr)_{\langle c\rangle\in\mathcal P_w}
\]
is an $r$-$(n,w,\lambda(q-1)^{s-1})$ design.
Here, the support of each one-dimensional subspace is counted exactly once.
\end{corollary}

\begin{proof}
Fix $R\in\binom{[n]}r$ and $S\in\binom Rs$.
By Theorem~\ref{thm:projective-characterization}, there are $\lambda$ one-dimensional subspaces that contain $R$ in their support and have the specified coordinate ratios.
There are $(q-1)^{s-1}$ possible tuples of coordinate ratios, and these partition the set of all one-dimensional subspaces whose supports contain $R$.
Thus, the number of blocks containing $R$ is $\lambda(q-1)^{s-1}$.
\end{proof}

In Corollary~\ref{cor:projective-supports}, by counting the blocks containing a fixed subset, we obtain the following divisibility condition \cite{BJL1999}.

\begin{corollary}\label{cor:linear-divisibility}
Let $1\leq s\leq r\leq w\leq n$ and suppose that $C_w$ is an $(r,s)$-$(n,w,q,\lambda)$ design.
Then for any $0\leq j\leq r$,
\[
\binom{w-j}{r-j}\mathrel{\Big|}
\lambda(q-1)^{s-1}\binom{n-j}{r-j}.
\]
\end{corollary}

\begin{proof}
In Corollary~\ref{cor:projective-supports}, counting in two ways the pairs $(R, B)$ where $R$ is an $r$-element subset containing a fixed $j$-element subset and $B$ is a block containing $R$, we find that the number of blocks containing the $j$-element subset is
\[
\lambda(q-1)^{s-1}
\frac{\binom{n-j}{r-j}}{\binom{w-j}{r-j}}.
\]
Since this number must be an integer, the above divisibility condition holds.
\end{proof}

\begin{example}
Let $n\geq2$ and consider the linear code
\[
C:=\{(c_1,\ldots,c_n)\in\mathbb F_q^n\mid c_1+\cdots+c_n=0\}.
\]
Codewords of weight $2$ are restricted to those having two nonzero coordinates with values $a$ and $-a$ for $a\in\mathbb F_q^\times$, so $|C_2|=(q-1)\binom n2$.
Fixing $R=\{i,j\}$, $S=\{i\}$, and $\omega\in\mathbb F_q^\times$, there is exactly one codeword satisfying $c_i=\omega$ and $c_j=-\omega$.
Therefore, $C_2$ is a $(2,1)$-$(n,2,q,1)$ design, and in the multiset of Corollary~\ref{cor:projective-supports}, the multiplicity of each $2$-element subset is $1$.

On the other hand, if we specify the nonzero values $\alpha,\beta\in\mathbb F_q^\times$ for two coordinates, the number of such codewords is $1$ if $\alpha+\beta=0$ and $0$ otherwise.
If $q>2$, both cases $\beta=-\alpha$ and $\beta\neq-\alpha$ can occur, which implies that $C_2$ is not a $(2,2)$-design.
\end{example}

\begin{example}
Nozaki--Watanabe provided a $(2,1)$-$(5,3,4,1)$ design \cite[Figure~1]{NW2025}.
However, a design with these parameters cannot be realized as $C_3$ of a linear code over $\mathbb F_4$.
Setting $j=0$ in Corollary~\ref{cor:linear-divisibility} requires $\binom32\mid\binom52$, but $3\nmid10$.
\end{example}

%%%%%%%%%%%%%%%%%%%%%%%%%%%%%%%%%%
\subsection{Character sums and dual codes}

The condition that tuples of nonzero coordinate values appear with the same frequency can be determined by verifying that the sums over all non-trivial characters vanish.
A character of a finite abelian group $G$ is a group homomorphism from $G$ to the multiplicative group $\mathbb C\setminus\{0\}$, and the set of all such characters is denoted by $\widehat G$.
The constant function $1$ is called the trivial character.
We use the orthogonality relations for characters \cite[Chapter~5]{LidlNiederreiter1997} in the following.

\begin{lemma}\label{lem:grp-character}
Let $G$ be a finite abelian group and $f:G\to\mathbb C$.
Then $f$ is constant if and only if for any $\rho\in\widehat G\setminus\{1\}$,
\[
\sum_{g\in G}f(g)\rho(g)=0.
\]
\end{lemma}

\begin{proof}
The inversion formula derived from the orthogonality of characters is given by
\[
f(g)=\frac1{|G|}\sum_{\rho\in\widehat G}
\left(\sum_{h\in G}f(h)\rho(h)\right)\overline{\rho(g)}.
\]
If the sums corresponding to the non-trivial characters all vanish, only the term for the trivial character remains, meaning that $f(g)$ does not depend on $g$.
Conversely, if $f$ is constant, the conclusion follows from the fact that the sum of the values of a non-trivial character over the group is $0$.
\end{proof}

Let $\widehat{\mathbb F_q^\times}$ denote the set of all characters of the multiplicative group $\mathbb F_q^\times$.
For $S\subseteq[n]$, the characters of the direct product group $(\mathbb F_q^\times)^S$ are given in the form
\[
(\omega_i)_{i\in S}\longmapsto\prod_{i\in S}\rho_i(\omega_i),
\qquad \rho_i\in\widehat{\mathbb F_q^\times}\quad(i\in S).
\]
Let $\mathcal X(S)$ be the set of tuples of characters for which this product is invariant under nonzero scalar multiplication:
\[
\mathcal X(S):=
\left\{(\rho_i)_{i\in S}\in
\bigl(\widehat{\mathbb F_q^\times}\bigr)^S\ \middle|\
\prod_{i\in S}\rho_i(a)=1\quad(a\in\mathbb F_q^\times)
\right\}.
\]

For $S\subseteq R\subseteq[n]$ with $|S|=s$, $|R|=r$, and $r\leq w\leq n$, and for a tuple of multiplicative characters $\boldsymbol\rho=(\rho_i)_{i\in S}$, we define
\[
\Phi_w(R,S;\boldsymbol\rho):=
\sum_{\substack{c\in C_w\\R\subseteq\operatorname{supp}(c)}}
\prod_{i\in S}\rho_i(c_i).
\]
Since $C_w$ is invariant under nonzero scalar multiplication, replacing the summation index $c$ with $ac$ yields
\begin{equation}\label{eq:scaling}
\Phi_w(R,S;\boldsymbol\rho)
=\left(\prod_{i\in S}\rho_i(a)\right)
\Phi_w(R,S;\boldsymbol\rho)
\qquad(a\in\mathbb F_q^\times).
\end{equation}
Therefore, if $\boldsymbol\rho\notin\mathcal X(S)$, the coefficient on the right-hand side is not $1$ for some $a$, which implies $\Phi_w(R,S;\boldsymbol\rho)=0$.

\begin{theorem}\label{thm:character-criterion}
Let $1\leq s\leq r\leq w\leq n$.
Then $C_w$ is an $(r,s)$-design if and only if the following two conditions hold:
\begin{enumerate}
\item $m_R(C_w)$ does not depend on the choice of $R\in\binom{[n]}r$.
\item For any $S\subseteq R\subseteq[n]$ with $|S|=s$ and $|R|=r$, and any $\boldsymbol\rho\in\mathcal X(S)\setminus\{(1)_{i\in S}\}$, we have $\Phi_w(R,S;\boldsymbol\rho)=0$.
\end{enumerate}
\end{theorem}

\begin{proof}
Fix $R$ and $S$, and apply Lemma~\ref{lem:grp-character} to $f(\omega)=m_{R,S}(C_w,\omega)$.
The sum for each character is
\[
\sum_{\omega\in(\mathbb F_q^\times)^S}
m_{R,S}(C_w,\omega)\prod_{i\in S}\rho_i(\omega_i)
=\Phi_w(R,S;\boldsymbol\rho),
\]
which is $0$ by Equation~\eqref{eq:scaling} if $\boldsymbol\rho\notin\mathcal X(S)$.
Thus, condition 2 is equivalent to saying that $m_{R,S}(C_w,\omega)$ does not depend on $\omega$.
In this case, taking the sum over $\omega$ yields
\[
m_{R,S}(C_w,\omega)=\frac{m_R(C_w)}{(q-1)^s},
\]
so this value also does not depend on $S$.
Combined with condition 1, it implies that this value is independent of $R$ as well, which means that $C_w$ is an $(r,s)$-design.
Conversely, if $C_w$ is an $(r,s)$-design, condition 1 holds by Lemma~\ref{lem:size}, and condition 2 also holds by the above equivalence.
\end{proof}

\begin{corollary}\label{cor:s-one-characterization}
Let $1\leq r\leq w\leq n$.
Then $C_w$ is an $(r,1)$-design if and only if $m_R(C_w)$ does not depend on the choice of $R\in\binom{[n]}r$.
\end{corollary}

\begin{proof}
When $|S|=1$, we have $\mathcal X(S)=\{(1)_{i\in S}\}$, so condition 2 of Theorem~\ref{thm:character-criterion} is always satisfied.
\end{proof}

Next, we express $\Phi_w$ as a sum over the codewords of $C^\perp$.
Fix a non-trivial character $\chi$ of the additive group $(\mathbb F_q,+)$.
For $\rho\in\widehat{\mathbb F_q^\times}$ and $b\in\mathbb F_q$, let
\[
G_\rho(b):=\sum_{a\in\mathbb F_q^\times}\rho(a)\chi(ab).
\]
From the orthogonality relations for characters \cite[Chapter~5]{LidlNiederreiter1997}, we have
\begin{equation}\label{eq:gauss-zero}
G_1(b)=
\begin{cases}q-1&b=0,\\-1&b\neq0,\end{cases}
\qquad
G_\rho(0)=0\quad(\rho\neq1).
\end{equation}

Let $N\geq0$ and $0\leq\ell\leq N$ be integers.
We set $\binom{x}{h}:=x(x-1)\cdots(x-h+1)/h!$ and $\binom{x}{0}:=1$, and define the $q$-ary Krawtchouk polynomial by
\[
P_\ell(x;N):=\sum_{h=0}^{\ell}(-1)^h(q-1)^{\ell-h}
\binom{x}{h}\binom{N-x}{\ell-h}
\]
\cite{MS1977}.
The generating function for integers $0\leq x\leq N$ is
\begin{equation}\label{eq:krawtchouk-generating}
\sum_{\ell=0}^N P_\ell(x;N)z^\ell
=(1+(q-1)z)^{N-x}(1-z)^x.
\end{equation}
Furthermore, the coefficient of $x^\ell$ is
\[
(-1)^\ell\sum_{h=0}^{\ell}
\frac{(q-1)^{\ell-h}}{h!(\ell-h)!}
=\frac{(-q)^\ell}{\ell!},
\]
so the degree of $P_\ell(x;N)$ is $\ell$.
We will use this property for interpolation in the next subsection.

Since $C=(C^\perp)^\perp$, the following equation follows from the orthogonality of additive characters \cite{MS1977}.

\begin{lemma}\label{lem:dual-character}
For any $v\in\mathbb F_q^n$,
\[
\frac1{|C^\perp|}\sum_{u\in C^\perp}\chi(u\cdot v)
=\begin{cases}1&v\in C,\\0&v\notin C.
\end{cases}
\]
\end{lemma}

\begin{proof}
If $v\in C$, then $u\cdot v=0$, so each term is $1$.
If $v\notin C=(C^\perp)^\perp$, then $u\cdot v\neq0$ for some $u\in C^\perp$.
As $a$ ranges over $\mathbb F_q$, we have $au\in C^\perp$, and $(au)\cdot v$ ranges over the entire $\mathbb F_q$.
Since $\chi$ is non-trivial, the mapping $u\mapsto\chi(u\cdot v)$ is a non-trivial character of the additive group $C^\perp$, and its sum vanishes.
\end{proof}

For $S\subseteq R\subseteq[n]$ with $|S|=s$ and $|R|=r$, and a tuple of multiplicative characters $\boldsymbol\rho=(\rho_i)_{i\in S}$, let
\[
H_j(R,S;\boldsymbol\rho):=
\sum_{\substack{u\in C^\perp\\
\operatorname{wt}(u|_{[n]\setminus R})=j}}
\left(\prod_{i\in S}G_{\rho_i}(u_i)\right)
\left(\prod_{i\in R\setminus S}G_1(u_i)\right)
\qquad(0\leq j\leq n-r).
\]
Similarly to the proof of the MacWilliams identity, we establish the following equation using the orthogonality of additive characters \cite{MS1977}.
For the Jacobi polynomials of codes, we refer the reader to Ozeki \cite{Ozeki1997}.
Also, for a MacWilliams-type transformation when variables are separated by coordinate sets and values, see \cite[Theorem~4.1]{CIT2024}.

\begin{proposition}\label{prop:dual-krawtchouk}
Let $1\leq s\leq r\leq w\leq n$, $S\subseteq R\subseteq[n]$ with $|S|=s$ and $|R|=r$, and let $\boldsymbol\rho=(\rho_i)_{i\in S}$ be a tuple of multiplicative characters.
Then
\[
\Phi_w(R,S;\boldsymbol\rho)
=\frac1{|C^\perp|}\sum_{j=0}^{n-r}
H_j(R,S;\boldsymbol\rho)P_{w-r}(j;n-r).
\]
\end{proposition}

\begin{proof}
Write $[z^k]F(z)$ for the coefficient of $z^k$ in $F(z)$.
By Lemma~\ref{lem:dual-character},
\[
\begin{aligned}
\Phi_w(R,S;\boldsymbol\rho)
&=\frac{1}{|C^\perp|}\sum_{u\in C^\perp}
  \sum_{\substack{v\in\mathbb F_q^n,\ \operatorname{wt}(v)=w\\
                  R\subseteq\operatorname{supp}(v)}}
  \chi(u\cdot v)\prod_{i\in S}\rho_i(v_i)\\
&=\frac{1}{|C^\perp|}\sum_{u\in C^\perp}
  \left(\prod_{i\in S}G_{\rho_i}(u_i)\right)
  \left(\prod_{i\in R\setminus S}G_1(u_i)\right)\\
&\qquad{}\times [z^{w-r}]
  \prod_{i\in[n]\setminus R}\bigl(1+G_1(u_i)z\bigr).
\end{aligned}
\]
For $i\in[n]\setminus R$, the choice $v_i=0$ contributes $1$,
and summing over $v_i\neq0$ contributes $G_1(u_i)z$.
Since all coordinates of $v$ on $R$ are nonzero,
the coefficient of $z^{w-r}$ selects the vectors of weight $w$.

Put $j=\operatorname{wt}(u|_{[n]\setminus R})$.
By Equations~\eqref{eq:gauss-zero} and~\eqref{eq:krawtchouk-generating},
\[
\begin{aligned}
[z^{w-r}]\prod_{i\in[n]\setminus R}\bigl(1+G_1(u_i)z\bigr)
&=[z^{w-r}](1+(q-1)z)^{n-r-j}(1-z)^j\\
&=P_{w-r}(j;n-r).
\end{aligned}
\]
Grouping the terms by $j$ gives the stated formula by the definition
of $H_j(R,S;\boldsymbol\rho)$.
\end{proof}

%%%%%%%%%%%%%%%%%%%%%%%%%%%%%%%%%%
\subsection{Sufficient conditions of Assmus--Mattson type}

The left-hand side of Proposition~\ref{prop:dual-krawtchouk} is $0$ at weights not occurring in $C$, providing linear equations for $H_j$.
First, we show by polynomial interpolation that all $H_j$ are $0$ if the number of nonzero $H_j$ is at most $d-r$.
Next, we also use equations derived from $C_w$, which is already known to be a design.
In what follows, the rank of a matrix is taken over $\mathbb C$.

\begin{theorem}\label{thm:abstract-AM}
Let $1\leq s\leq r<d$.
Assume that for any $S\subseteq R\subseteq[n]$ with $|S|=s$ and $|R|=r$, and for any $\boldsymbol\rho\in\mathcal X(S)\setminus\{(1)_{i\in S}\}$, we have
\[
\bigl|\{j\in\{0,\ldots,n-r\}\mid H_j(R,S;\boldsymbol\rho)\neq0\}\bigr|\leq d-r.
\]
Then, for these $R,S,\boldsymbol\rho$,
\[
\Phi_w(R,S;\boldsymbol\rho)=0\qquad(r\leq w\leq n)
\]
Therefore, for any $w$ with $r\leq w\leq n$ such that $m_R(C_w)$ is independent of $R\in\binom{[n]}r$, $C_w$ is an $(r,s)$-design.
\end{theorem}

\begin{proof}
If $s=1$, then $\mathcal X(S)\setminus\{(1)_{i\in S}\}=\varnothing$, and when $m_R(C_w)$ is independent of $R$, Corollary~\ref{cor:s-one-characterization} implies that $C_w$ is an $(r,1)$-design.
Suppose $s\geq2$ and fix $R,S,\boldsymbol\rho$.
Let
\[
I:=\{j\in\{0,\ldots,n-r\}\mid H_j(R,S;\boldsymbol\rho)\neq0\},
\qquad h:=|I|.
\]
If $h=0$, then Proposition~\ref{prop:dual-krawtchouk} implies that all $\Phi_w(R,S;\boldsymbol\rho)$ are $0$.
In the following, we assume $h\geq1$.
Since $h\leq d-r$, for $0\leq\ell<h$ we have $1\leq r+\ell<d$, and thus $C_{r+\ell}=\varnothing$.
Therefore, Proposition~\ref{prop:dual-krawtchouk} yields
\[
\sum_{j\in I}H_j(R,S;\boldsymbol\rho)P_\ell(j;n-r)=0
\qquad(0\leq\ell<h).
\]
Since the degrees of $P_0(x;n-r),\ldots,P_{h-1}(x;n-r)$ are $0,\ldots,h-1$, respectively, they form a basis for the vector space of polynomials of degree at most $h-1$.
Thus, for any such polynomial $f$,
\[
\sum_{j\in I}H_j(R,S;\boldsymbol\rho)f(j)=0
\]
holds.
For each $j\in I$, using the interpolation polynomial
\[
f(x)=\prod_{k\in I\setminus\{j\}}\frac{x-k}{j-k},
\]
we obtain $H_j(R,S;\boldsymbol\rho)=0$, which contradicts $j\in I$.
Therefore, $I=\varnothing$.
By Proposition~\ref{prop:dual-krawtchouk}, $\Phi_w(R,S;\boldsymbol\rho)=0$ for all $r\leq w\leq n$.
Consequently, when $m_R(C_w)$ is independent of $R$, Theorem~\ref{thm:character-criterion} ensures that $C_w$ is an $(r,s)$-design.
\end{proof}

For $C_w$ to be an $(r,s)$-design, in addition to $\Phi_w(R,S;\boldsymbol\rho)=0$ for $\boldsymbol\rho\in\mathcal X(S)\setminus\{(1)_{i\in S}\}$, $m_R(C_w)$ must also be independent of $R$.
To show that $m_R(C_w)$ is independent of $R$, we use the Assmus--Mattson theorem \cite[Theorem~1.1]{MT2018}, which counts the support sets, counting each codeword exactly once.
Let the set of nonzero weights occurring in the dual code be
\[
\mathcal W^\perp:=\{\operatorname{wt}(u)\mid u\in C^\perp\setminus\{0\}\}.
\]

\begin{lemma}\label{lem:ordinary-AM}
Let $1\leq r<d$, and assume
\[
|\mathcal W^\perp\cap\{1,\ldots,n-r\}|\leq d-r.
\]
Then, for any $d\leq w\leq n$, $m_R(C_w)$ is independent of $R\in\binom{[n]}r$.
\end{lemma}

\begin{proof}
If $C_w\neq\varnothing$, then by \cite[Theorem~1.1]{MT2018}, the support sets of the codewords, counted with multiplicity, form an $r$-design.
Since the number of blocks containing $R$ is $m_R(C_w)$, we obtain the conclusion.
If $C_w=\varnothing$, then $m_R(C_w)=0$.
\end{proof}

Next, we define a set containing indices for which $H_j\neq0$ from the weights occurring in the dual code.
For $1\leq r\leq n$, let
\[
J:=\bigcup_{\ell\in\mathcal W^\perp}
\{j\in\{0,\ldots,n-r\}\mid \ell-r\leq j\leq\ell-2\}.
\]
If $\boldsymbol\rho\in\mathcal X(S)\setminus\{(1)_{i\in S}\}$, at least two of the $\rho_i$ are nontrivial.
From Equation~\eqref{eq:gauss-zero}, $u\in C^\perp$ corresponding to a nonzero term of $H_j(R,S;\boldsymbol\rho)$ is nonzero at those coordinates.
Therefore, at least two components are nonzero on $R$, and $j\leq\operatorname{wt}(u)-2$.

\begin{lemma}\label{lem:gauss-support}
Let $1\leq s\leq r\leq n$, $S\subseteq R\subseteq[n]$, $|S|=s$, $|R|=r$, and let $\boldsymbol\rho\in\mathcal X(S)\setminus\{(1)_{i\in S}\}$.
If $H_j(R,S;\boldsymbol\rho)\neq0$, then $j\in J$.
\end{lemma}

\begin{proof}
Let $T:=\{i\in S\mid\rho_i\neq1\}$.
Since $\boldsymbol\rho$ is not the tuple of trivial characters and its product is the trivial character, we have $|T|\geq2$.
From Equation~\eqref{eq:gauss-zero}, $u\in C^\perp$ corresponding to a nonzero term of $H_j(R,S;\boldsymbol\rho)$ satisfies $T\subseteq\operatorname{supp}(u)$.
If $H_j(R,S;\boldsymbol\rho)\neq0$, there exists at least one such term, and the corresponding $u$ satisfies
\[
\operatorname{wt}(u|_{[n]\setminus R})=j,
\qquad j+2\leq\operatorname{wt}(u)\leq j+r.
\]
Thus, $\operatorname{wt}(u)\in\mathcal W^\perp$ and $j\in J$.
\end{proof}

If $C_w=\varnothing$ or $C_w$ is an $(r,s)$-design, then $\Phi_w(R,S;\boldsymbol\rho)=0$ for $\boldsymbol\rho\in\mathcal X(S)\setminus\{(1)_{i\in S}\}$.
By Proposition~\ref{prop:dual-krawtchouk}, we obtain linear equations for $H_j$ from these equalities.

\begin{theorem}\label{thm:rank}
Let $1\leq s\leq r<d$, and assume
\[
|\mathcal W^\perp\cap\{1,\ldots,n-r\}|\leq d-r.
\]
Suppose that for each $m\in L\subseteq\{0,\ldots,n-r\}$, $C_{r+m}=\varnothing$ or $C_{r+m}$ is an $(r,s)$-design.
Furthermore, assume
\[
\operatorname{rank}\bigl(P_m(j;n-r)\bigr)_{m\in L,\ j\in J}=|J|.
\]
Then, for any $d\leq w\leq n$, $C_w$ is an $(r,s)$-design.
\end{theorem}

\begin{proof}
The case $s=1$ follows from Lemma~\ref{lem:ordinary-AM} and Corollary~\ref{cor:s-one-characterization}.
Suppose $s\geq2$, and fix $S\subseteq R\subseteq[n]$ with $|S|=s$, $|R|=r$, and $\boldsymbol\rho\in\mathcal X(S)\setminus\{(1)_{i\in S}\}$.
For each $m\in L$, we have $\Phi_{r+m}(R,S;\boldsymbol\rho)=0$ by definition if $C_{r+m}=\varnothing$, and by Theorem~\ref{thm:character-criterion} if $C_{r+m}$ is a design.
Therefore, Proposition~\ref{prop:dual-krawtchouk} and Lemma~\ref{lem:gauss-support} yield
\[
\sum_{j\in J}H_j(R,S;\boldsymbol\rho)P_m(j;n-r)=0
\qquad(m\in L).
\]
Since the coefficient matrix has rank equal to the number of columns, $H_j(R,S;\boldsymbol\rho)=0$ for $j\in J$.
For $j\in\{0,\ldots,n-r\}\setminus J$, we also have $H_j(R,S;\boldsymbol\rho)=0$ by Lemma~\ref{lem:gauss-support}.
Thus, Proposition~\ref{prop:dual-krawtchouk} implies that $\Phi_w(R,S;\boldsymbol\rho)=0$ for all $r\leq w\leq n$.
On the other hand, Lemma~\ref{lem:ordinary-AM} ensures that $m_R(C_w)$ is independent of $R$.
Therefore, both conditions of Theorem~\ref{thm:character-criterion} are satisfied.
\end{proof}

\begin{proof}
Apply Theorem~\ref{thm:rank} for each $s$.
\end{proof}

These conditions are determined by $q,n,r$ and the weights occurring in $C$ and $C^\perp$.
By Theorem~\ref{thm:abstract-AM}, we may assume $|J|\leq d-r$ instead of the rank condition.

\begin{corollary}\label{cor:ordinary-weight}
Let $1\leq r<d$.
If
\[
|\mathcal W^\perp\cap\{1,\ldots,n-r\}|\leq d-r,
\qquad |J|\leq d-r,
\]
then for any $1\leq s\leq r$ and $d\leq w\leq n$, $C_w$ is an $(r,s)$-design, and its index is
\[
\frac{|C_w|\binom wr}{(q-1)^s\binom nr}.
\]
\end{corollary}

\begin{proof}
Fix $1\leq s\leq r$ and $d\leq w\leq n$.
The first hypothesis and Lemma~\ref{lem:ordinary-AM} imply that
$m_R(C_w)$ is independent of $R$.

For each $S\subseteq R\subseteq[n]$ with $|S|=s$ and $|R|=r$,
and each $\boldsymbol\rho\in\mathcal X(S)\setminus\{(1)_{i\in S}\}$,
Lemma~\ref{lem:gauss-support} gives
\[
\bigl\{j\in\{0,\ldots,n-r\}\mid
       H_j(R,S;\boldsymbol\rho)\neq0\bigr\}\subseteq J.
\]
Since $|J|\leq d-r$, Theorem~\ref{thm:abstract-AM} implies that
$C_w$ is an $(r,s)$-design. Its index is given by Lemma~\ref{lem:size}.
\end{proof}

\begin{corollary}\label{cor:few-dual-weights}
Let $2\leq r<d$.
If $(r-1)|\mathcal W^\perp|\leq d-r$, then for any $1\leq s\leq r$ and $d\leq w\leq n$, $C_w$ is an $(r,s)$-design.
In particular, for $r=2$, $|\mathcal W^\perp|\leq d-2$ is sufficient.
\end{corollary}

\begin{proof}
For each $\ell\in\mathcal W^\perp$, there are at most $r-1$ integers satisfying $\ell-r\leq j\leq\ell-2$.
Therefore,
\[
|J|\leq(r-1)|\mathcal W^\perp|,
\qquad
|\mathcal W^\perp\cap\{1,\ldots,n-r\}|\leq|\mathcal W^\perp|.
\]
Since $r-1\geq1$, both conditions of Corollary~\ref{cor:ordinary-weight} are satisfied.
\end{proof}

\begin{example}\label{ex:hexacode}
Let $\mathbb F_4=\{0,1,\omega,\omega^2\}$ with $\omega^2+\omega+1=0$.
The code $C=\{uG\mid u\in\mathbb F_4^3\}$ with the generator matrix
\[
G=\begin{pmatrix}
1&0&0&1&\omega&\omega\\
0&1&0&\omega&1&\omega\\
0&0&1&\omega&\omega&1
\end{pmatrix}
\]
is the Hexacode \cite{MOWW1978}.
By direct calculation, both $C$ and $C^\perp$ have $1$, $45$, and $18$ codewords of weights $0$, $4$, and $6$, respectively, and no other weights occur.
Thus, $d=4$, and for $r=3$, we have
\[
\mathcal W^\perp=\{4,6\},
\qquad |\mathcal W^\perp\cap\{1,2,3\}|=0,
\qquad J=\{1,2,3\}.
\]

First, we show that $C_6$ is a $(3,2)$-design.
The representatives of the six one-dimensional subspaces of $\mathcal P_6$ with the first component being $1$ are the rows of the following matrix:
\[
\begin{pmatrix}
1&1&1&1&1&1\\
1&1&\omega^2&\omega&\omega&\omega^2\\
1&\omega&\omega&1&\omega^2&\omega^2\\
1&\omega&\omega^2&\omega^2&1&\omega\\
1&\omega^2&1&\omega&\omega^2&\omega\\
1&\omega^2&\omega&\omega^2&\omega&1
\end{pmatrix}.
\]
For any $1\leq a<b\leq6$, the ratios $c_b/c_a$ of these six representatives contain $1,\omega,\omega^2$ twice each.
Since all support sets are $[6]$, Theorem~\ref{thm:projective-characterization} implies that $C_6$ is a $(3,2)$-$(6,6,4,2)$ design.

We apply Theorem~\ref{thm:rank} for $r=3$, $s=2$, and $L=\{0,2,3\}$.
We have $C_3=C_5=\varnothing$, and $C_6$ is a $(3,2)$-design.
Furthermore, for $q=4$, the determinant of
\[
\bigl(P_m(j;3)\bigr)_{m\in\{0,2,3\},\ j\in\{1,2,3\}}
=\begin{pmatrix}1&1&1\\3&-5&3\\-9&3&-1\end{pmatrix}
\]
is $-64$, so the rank condition is also satisfied.
Therefore, $C_4$ is a $(3,2)$-design, and its index by Lemma~\ref{lem:size} is
\[
\frac{|C_4|\binom43}{(4-1)^2\binom63}
=\frac{45\cdot4}{9\cdot20}=1.
\]
Thus, $C_4$ is a $(3,2)$-$(6,4,4,1)$ design.
\end{example}

\begin{remark}\label{rem:new-few-dual-weights}
Let $2\leq r<d$.
The condition of Corollary~\ref{cor:few-dual-weights} is never satisfied for $r\geq3$.
In fact, Delsarte's bound \cite[Section~2.1]{BZ2025} gives
\[
\max_{x\in\mathbb F_q^n}\min_{c\in C}d_H(x,c)\leq|\mathcal W^\perp|.
\]
On the other hand, for a vector $x$ of weight $\lfloor d/2\rfloor$, the distance to the zero codeword is $\lfloor d/2\rfloor$, and for any nonzero codeword $c$,
\[
d_H(x,c)\geq\operatorname{wt}(c)-\operatorname{wt}(x)
\geq d-\lfloor d/2\rfloor\geq\lfloor d/2\rfloor.
\]
Therefore, $\lfloor d/2\rfloor\leq|\mathcal W^\perp|$, and if $r\geq3$, then
\[
(r-1)|\mathcal W^\perp|\geq2\lfloor d/2\rfloor\geq d-1>d-r.
\]
The conclusion for $r=2$ also follows from Wu--Zhou \cite[Theorem~3.6]{WZ2026}.
\end{remark}

%%%%%%%%%%%%%%%%%%%%%%%%%%%%%%%%%%
\subsection{Necessary and sufficient conditions for two-weight codes}\label{subsec:new-two-weight}

In this subsection, let $q>2$, $0\neq C\subsetneq\mathbb F_q^n$, and $D:=C^\perp$.
We denote the set of codewords of weight $v$ in $D$ by $D_v$, and its minimum distance by $\delta$:
\[
\begin{aligned}
D_v&:=\{u\in D\mid\operatorname{wt}(u)=v\}\quad(1\leq v\leq n),\\
\delta&:=d(C^\perp)=\min\{\operatorname{wt}(u)\mid u\in D\setminus\{0\}\}.
\end{aligned}
\]
Both $C$ and $D$ are nonzero, and $d$ continues to denote the minimum distance of $C$.

In a $(t,t)$-design, nonzero values are fixed at all $t$ specified coordinates.

\begin{definition}
Let $0\leq t\leq w\leq n$ and $\lambda\geq0$ be integers.
If $Y\subseteq W_w(n,q)$ is a $(t,t)$-$(n,w,q,\lambda)$ design, we call $Y$ a $q$-ary $t$-$(n,w,\lambda)$ design.
That is, the condition is
\[
|\{y\in Y\mid y|_I=a\}|=\lambda
\qquad\left(I\in\binom{[n]}t,\ a\in(\mathbb F_q^\times)^I\right).
\]
\end{definition}

A linear code $C$ satisfying
\[
\bigl|\{\operatorname{wt}(c)\mid c\in C\setminus\{0\}\}\bigr|=2
\]
is called a two-weight code.
For two-weight codes, we show that if a nonempty $C_w$ is a $q$-ary $t$-design with $t\geq2$, then $\delta\geq t+2$.
We first present two lemmas used in the proof.
The index when reducing the strength is given by the following formula \cite[Lemma~2.3]{WZ2026}.

\begin{lemma}\label{lem:new-descent}
Let $0\leq t\leq w\leq n$, and suppose $Y\subseteq W_w(n,q)$ is a $q$-ary $t$-$(n,w,\lambda)$ design.
Then for any $0\leq j\leq t$, $Y$ is a $q$-ary $j$-design, and its index is
\[
\lambda_j=(q-1)^{t-j}
\frac{\binom{n-j}{t-j}}{\binom{w-j}{t-j}}\lambda.
\]
\end{lemma}

\begin{proof}
Fix $J\in\binom{[n]}j$ and $a\in(\mathbb F_q^\times)^J$.
For each $t$-element subset containing $J$, summing the number of vectors over all possible assignments of nonzero values to the coordinates outside $J$ gives
\[
\binom{w-j}{t-j}|\{y\in Y\mid y|_J=a\}|
=\binom{n-j}{t-j}(q-1)^{t-j}\lambda.
\]
On the left-hand side, each $y$ is counted with a multiplicity equal to the number of $t$-element subsets containing $J$ and contained in $\operatorname{supp}(y)$.
Dividing by $\binom{w-j}{t-j}>0$, we obtain $|\{y\in Y\mid y|_J=a\}|=\lambda_j$.
\end{proof}

In the next lemma, we specify values (including zeros) on $t$ coordinates and apply a MacWilliams-type transform to a polynomial representing the weight of the remaining components as exponents.

\begin{lemma}\label{lem:new-qary-duality}
Let $1\leq t\leq n$, $d\geq t$, and $\delta\geq t$.
All nonempty $C_w$ for $1\leq w\leq n$ are $q$-ary $t$-designs if and only if all nonempty $D_v$ for $1\leq v\leq n$ are $q$-ary $t$-designs.
\end{lemma}

\begin{proof}
Assume that all nonempty $C_w$ are $q$-ary $t$-designs.
We will show that all nonempty $D_v$ are also $q$-ary $t$-designs.
For $I\in\binom{[n]}t$ and $a\in\mathbb F_q^I$, using indeterminates $x,y$, let
\[
L^C_I(a;x,y):=
\sum_{\substack{c\in C\\c|_I=a}}
x^{n-t-\operatorname{wt}(c|_{[n]\setminus I})}
y^{\operatorname{wt}(c|_{[n]\setminus I})}.
\]
We define $L^D_I$ by replacing $C$ with $D$.
First, we show that $L^C_I(a;x,y)$ is independent of $I$ and depends only on the weight of $a$.

Fix a nonempty $C_w$ with $w>0$.
Since $w\geq d\geq t$, Lemma~\ref{lem:new-descent} implies that the number of codewords with specified nonzero values on $j$ coordinates is independent of the choice of coordinates and values.
We denote this number by $\mu_j$ for $0\leq j\leq t$.
By the principle of inclusion-exclusion, the number of codewords with $j$ specified nonzero values and $t-j$ specified zeros on $I$ is
\[
\sum_{h=0}^{t-j}(-1)^h\binom{t-j}{h}(q-1)^h\mu_{j+h}.
\]
This is obtained by selecting $h$ coordinates from those specified as zero, varying their nonzero values, and summing the number of codewords satisfying the specification.
This number, as well as the weight $w-j$ of the remaining coordinates, depends only on $j$.
Since the zero codeword is included in the sum only when $a=0$, $L^C_I(a;x,y)$ depends only on $j=\operatorname{wt}(a)$.
We write this common polynomial as $F_j(x,y)$.

Using a fixed additive character $\chi$, let $X:=x+(q-1)y$ and $Y:=x-y$.
Applying Lemma~\ref{lem:dual-character} to $D$ and using $D^\perp=C$, we have
\[
\begin{aligned}
L^D_I(b;x,y)
&=\frac1{|C|}\sum_{a\in\mathbb F_q^I}
\chi\!\left(\sum_{i\in I}a_ib_i\right)L^C_I(a;X,Y)\\
&=\frac1{|C|}\sum_{j=0}^tP_j(\operatorname{wt}(b);t)F_j(X,Y)
\qquad(b\in\mathbb F_q^I).
\end{aligned}
\]
In the first equality, the sum is taken over each component outside $I$.
For each $c\in C$, the character sum over components $i\notin I$ is $X$ if $c_i=0$ and $Y$ if $c_i\neq0$.
In the second equality, terms are grouped by the weight of $a$.
From Equation~\eqref{eq:krawtchouk-generating}, the character sum with respect to $a$ of weight $j$ is equal to $P_j(\operatorname{wt}(b);t)$.

Therefore, $L^D_I(b;x,y)$ is independent of $I$ and depends only on the weight of $b$.
Let $v\geq\delta\geq t$ and take $b\in(\mathbb F_q^\times)^I$.
The coefficient of $x^{n-v}y^{v-t}$ is $|\{u\in D_v\mid u|_I=b\}|$, which is independent of $I$ and $b$.
Thus, $D_v$ is a $q$-ary $t$-design.
Since $D^\perp=C$ and the minimum distances of both codes are at least $t$, the converse can be proved similarly by interchanging $C$ and $D$.
\end{proof}

For two-weight codes, $\delta\geq t+2$ is a sufficient condition to yield a $q$-ary $t$-design \cite[Theorem~3.6]{WZ2026}.
The following theorem establishes the necessity of this condition from the fact that a single nonempty $C_w$ is a $q$-ary $t$-design.

\begin{theorem}\label{thm:new-two-weight-exact}
Let $q>2$ be a prime power, $0\neq C\subsetneq\mathbb F_q^n$ be a two-weight linear code, and $t\geq2$ be an integer.
Then we have
\[
\begin{gathered}
\exists\,w\in\{1,\ldots,n\},\ w\geq t:\quad
C_w\neq\varnothing,\quad C_w\text{ is }q\text{-ary }t\text{-design}
\\
\Longleftrightarrow\quad d(C^\perp)\geq t+2.
\end{gathered}
\]
When these conditions hold, $d>t$, and all nonempty $C_w$ and $D_v$ for $1\leq w,v\leq n$ are $q$-ary $t$-designs.
\end{theorem}

\begin{proof}
If $t>n$, the left-hand side is false, and the right-hand side is also false since $\delta\leq n<t+2$.
In the following, assume $2\leq t\leq n$.
To prove necessity, let the two nonzero weights of $C$ be $w,h$, and assume $w\geq t$ and that $C_w$ is a $q$-ary $t$-design.
Since $C_w\neq\varnothing$, its index $\lambda$ is positive by Lemma~\ref{lem:size}.
Let $k:=\dim_{\mathbb F_q}C$.

First, we show that for any $I\in\binom{[n]}t$, the linear map
\[
C\longrightarrow\mathbb F_q^I,\qquad c\longmapsto c|_I
\]
is surjective.
By the design condition and $\lambda>0$, its image contains $(\mathbb F_q^\times)^I$.
Since $q>2$, we can choose two vectors from $(\mathbb F_q^\times)^I$ that differ in exactly one coordinate's nonzero value.
Their difference is nonzero only at that coordinate.
By taking this difference for each coordinate, the image contains vectors that span $\mathbb F_q^I$, making the map surjective.
Thus, the preimage of each element has cardinality $q^{k-t}$.
If there exists $0\neq u\in D$ with $\operatorname{wt}(u)\leq t$, we can choose $I$ such that $\operatorname{supp}(u)\subseteq I$.
In this case, the image of the map satisfies the nontrivial linear relation $\sum_{i\in I}u_ic_i=0$, contradicting surjectivity.
Therefore, $\delta>t$.

Next, we show $d>t$.
If $w=t$, specifying the nonzero values on $t$ coordinates uniquely determines a vector of weight $t$.
Since $\lambda>0$, all vectors of weight $t$ belong to $C$.
By choosing a support set containing a given coordinate and taking the difference of two vectors whose nonzero value differs only at that coordinate, $C$ contains vectors that are nonzero only at each individual coordinate.
This implies $C=\mathbb F_q^n$, which contradicts the assumption.
Therefore, $w>t$.

Assume $h\leq t$.
Since any $H\in\binom{[n]}h$ is contained in a $t$-element subset, the restriction map $c\mapsto c|_H$ onto $H$ is also surjective.
The number of codewords with specified nonzero values on $H$ is $q^{k-h}$ for the whole code $C$, and by Lemma~\ref{lem:new-descent}, this number in $C_w$ is independent of the choice of coordinates and values.
Since the zero codeword does not have the specified nonzero values and the only nonzero weights are $w,h$, their difference shows that the number of such codewords in $C_h$ is also independent of coordinates and values.
Because a codeword of weight $h$ is uniquely determined by specifying its nonzero values on $H$, this number is either $0$ or $1$.
Since $C_h\neq\varnothing$, the common value is $1$, meaning all vectors of weight $h$ belong to $C$.
Taking the difference of two vectors that differ only in one coordinate's nonzero value similarly yields $C=\mathbb F_q^n$.
Thus, $h>t$, which implies $d>t$.

Since the restriction map to any $t$ coordinates is surjective, the number of codewords in $C_h$ with specified nonzero values is always $q^{k-t}-\lambda$.
Therefore, $C_h$ is also a $q$-ary $t$-design.
By $d>t$, $\delta>t$, and Lemma~\ref{lem:new-qary-duality}, all nonempty $D_v$ are also $q$-ary $t$-designs.
We will now eliminate the remaining possibility $\delta=t+1$.

Assume $t\geq3$.
Take $u\in D_{t+1}$ and $S\in\binom{\operatorname{supp}(u)}t$, and change one nonzero value of $u|_S$ to another nonzero value.
Since $D_{t+1}$ is a nonempty $q$-ary $t$-design, its index is positive, so there exists $u'\in D_{t+1}$ taking the modified values.
$u$ and $u'$ differ in exactly one coordinate on $S$ and each has one nonzero component outside $S$.
Therefore,
\[
0<\operatorname{wt}(u-u')\leq3<t+1=\delta.
\]
Since $u-u'\in D$, this contradicts the definition of the minimum distance.

Assume $t=2$ and $\delta=3$.
Let the columns of the generator matrix of $C$ be $g_1,\ldots,g_n\in\mathbb F_q^k$.
A zero column would give a codeword of weight $1$ in $D$, and two proportional columns would give a codeword of weight $2$.
Thus, all columns are nonzero, and any two columns are linearly independent.
For $0\neq g\in\mathbb F_q^k$, let $\langle g\rangle:=\{ag\mid a\in\mathbb F_q\}$, and set
\[
P:=\{\langle g_i\rangle\mid i\in[n]\}.
\]

Fix distinct $i,j\in[n]$ and $\alpha,\beta\in\mathbb F_q^\times$.
For $u\in D_3$ satisfying $u_i=\alpha$ and $u_j=\beta$, the remaining nonzero coordinate $\ell$ and its value $\gamma$ satisfy
\[
\alpha g_i+\beta g_j+\gamma g_\ell=0.
\]
Since $g_i,g_j$ are linearly independent and $\alpha,\beta\neq0$, $\langle\alpha g_i+\beta g_j\rangle$ is distinct from $\langle g_i\rangle,\langle g_j\rangle$.
Because each one-dimensional subspace contains at most one column, $\ell,\gamma$ satisfying the above equation are unique if they exist.
Therefore,
\[
|\{u\in D_3\mid u_i=\alpha,\ u_j=\beta\}|
=\begin{cases}
1&\langle\alpha g_i+\beta g_j\rangle\in P,\\
0&\langle\alpha g_i+\beta g_j\rangle\notin P.
\end{cases}
\]
Since $D_3$ is a nonempty $q$-ary $2$-design, this number is always $1$.
As $\alpha,\beta$ vary, we obtain all one-dimensional subspaces contained in the two-dimensional subspace spanned by $g_i,g_j$, other than $\langle g_i\rangle,\langle g_j\rangle$.
Thus, all one-dimensional subspaces contained in the two-dimensional subspace spanned by any two distinct elements of $P$ must belong to $P$.

Let $V:=\bigcup_{L\in P}L$.
$V$ contains zero and is closed under scalar multiplication.
Since all one-dimensional subspaces contained in the two-dimensional subspace spanned by two linearly independent vectors in $V$ belong to $P$, their sum also belongs to $V$.
The sum of two linearly dependent vectors belongs to the same one-dimensional subspace or is zero.
Therefore, $V$ is a linear subspace.
Since the columns of the generator matrix span $\mathbb F_q^k$, we have $V=\mathbb F_q^k$.
This means the generator matrix has exactly one nonzero vector from each one-dimensional subspace of $\mathbb F_q^k$ as a column.

Each nonzero codeword is obtained by evaluating a nonzero linear functional on these columns.
Its kernel is $(k-1)$-dimensional and contains $(q^{k-1}-1)/(q-1)$ one-dimensional subspaces.
Therefore, the weight of each nonzero codeword is equal to
\[
\frac{q^k-1}{q-1}-\frac{q^{k-1}-1}{q-1}=q^{k-1}.
\]
This contradicts the assumption that it is a two-weight code.
Consequently, $\delta=t+1$ cannot occur, yielding $\delta\geq t+2$.

Finally, we prove sufficiency.
Assume $\delta\geq t+2$.
Since the dimension of $D$ is $n-k$, the Singleton bound \cite{MS1977} gives $\delta\leq k+1$, hence $k\geq\delta-1$.
The map to the coordinates corresponding to $k$ linearly independent columns of the generator matrix is surjective, so there exists a codeword where all these $k$ components are $1$.
Thus, the maximum nonzero weight of $C$ is at least $k$, and in particular at least $\delta-2$.
Since there are two nonzero weights and $\delta>2$, the condition of Wu--Zhou \cite[Theorem~3.6]{WZ2026} is satisfied.
The strength given by this theorem is $\max\{d-|\mathcal W^\perp|,\delta-2\}$.
This value is at least $\delta-2$ and at most the maximum nonzero weight of $C$.
Therefore, by Lemma~\ref{lem:new-descent}, the set of all codewords of maximum nonzero weight in $C$ is a $q$-ary $(\delta-2)$-design, and since $\delta-2\geq t$, it is also a $q$-ary $t$-design.
From the proof of necessity, $d>t$ and all nonempty $C_w,D_v$ are $q$-ary $t$-designs.
\end{proof}

\begin{corollary}\label{cor:new-WZ-open-problem}
Let $q>2$ be a prime power, and $0\neq C\subsetneq\mathbb F_q^n$ be a two-weight linear code with $d(C^\perp)=3$.
Then, for any integers $w,t$ satisfying $C_w\neq\varnothing$ and $2\leq t\leq w\leq n$, $C_w$ is not a $q$-ary $t$-design.
\end{corollary}

\begin{proof}
If such a $C_w$ exists, Theorem~\ref{thm:new-two-weight-exact} implies $d(C^\perp)\geq t+2\geq4$, which contradicts the assumption.
\end{proof}

This corollary provides a negative answer to the open problem posed by Wu--Zhou \cite[p.~16, after Remark~4.15]{WZ2026}, asking whether a two-weight linear code with dual minimum distance $3$ can have a nonempty $C_w$ that forms a $q$-ary $t$-design for $t\geq2$.

%%%%%%%%%%%%%%%%%%%%%%%%%%%%%%
\section*{Acknowledgements}
The research of Ryutaro Misawa was supported by JSPS Grant-in-Aid for JSPS Fellows (26KJ0570).
The research of Yuuho Tanaka was supported by JSPS KAKENHI Grant Number JP25K23339.

%%%%%%%%%%%%%%%%%%%%%%%%%%%%%%


\begin{thebibliography}{99}

\bibitem{AM1969}
E. F. Assmus, Jr. and H. F. Mattson, Jr.,
\newblock New $5$-designs,
\newblock \emph{Journal of Combinatorial Theory} \textbf{6} (1969), 122--151.
\newblock DOI: \url{https://doi.org/10.1016/S0021-9800(69)80115-8}.

\bibitem{Bahmanian2014}
A. Bahmanian,
\newblock Connected Baranyai's theorem,
\newblock \emph{Combinatorica} \textbf{34} (2014), 129--138.
\newblock DOI: \url{https://doi.org/10.1007/s00493-014-2928-3}.

\bibitem{Baranyai1975}
Z. Baranyai,
\newblock On the factorization of the complete uniform hypergraph,
\newblock in \emph{Infinite and Finite Sets}, Colloquia Mathematica Societatis J\'anos Bolyai, Vol.~10, North-Holland, Amsterdam, 1975, 91--108.

\bibitem{BJL1999}
T. Beth, D. Jungnickel, and H. Lenz,
\newblock \emph{Design Theory, Vol. I},
\newblock second edition, Cambridge University Press, Cambridge, 1999.

\bibitem{BZ2025}
J. Borges and V. Zinoviev,
\newblock On self-dual completely regular codes with covering radius $\rho\leq3$,
\newblock \emph{Finite Fields and Their Applications} \textbf{105} (2025), Article 102617.
\newblock DOI: \url{https://doi.org/10.1016/j.ffa.2025.102617}.

\bibitem{CIT2024}
H. S. Chakraborty, R. Ishikawa, and Y. Tanaka,
\newblock Jacobi polynomials and design theory II,
\newblock \emph{Discrete Mathematics} \textbf{347} (2024), Article 113818.
\newblock DOI: \url{https://doi.org/10.1016/j.disc.2023.113818}.

\bibitem{CheeLing2007}
Y. M. Chee and S. Ling,
\newblock Constructions for $q$-ary constant-weight codes,
\newblock \emph{IEEE Transactions on Information Theory} \textbf{53} (2007), no.~1, 135--146.
\newblock DOI: \url{https://doi.org/10.1109/TIT.2006.887499}.

\bibitem{Delsarte1973}
P. Delsarte,
\newblock \emph{An Algebraic Approach to the Association Schemes of Coding Theory},
\newblock Philips Research Reports Supplements, no.~10 (1973), 1--97.

\bibitem{Etzion2022}
T. Etzion,
\newblock Non-binary diameter perfect constant-weight codes,
\newblock \emph{IEEE Transactions on Information Theory} \textbf{68} (2022), no.~2, 891--904.
\newblock DOI: \url{https://doi.org/10.1109/TIT.2021.3125277}.

\bibitem{HSS1999}
A. S. Hedayat, N. J. A. Sloane, and J. Stufken,
\newblock \emph{Orthogonal Arrays: Theory and Applications},
\newblock Springer, New York, 1999.

\bibitem{Ji2005}
L. Ji,
\newblock A new existence proof for large sets of disjoint Steiner triple systems,
\newblock \emph{Journal of Combinatorial Theory, Series A} \textbf{112} (2005), 308--327.
\newblock DOI: \url{https://doi.org/10.1016/j.jcta.2005.06.005}.

\bibitem{LidlNiederreiter1997}
R. Lidl and H. Niederreiter,
\newblock \emph{Finite Fields},
\newblock second edition, Cambridge University Press, Cambridge, 1997.

\bibitem{MOWW1978}
F. J. MacWilliams, A. M. Odlyzko, N. J. A. Sloane, and H. N. Ward,
\newblock Self-dual codes over GF(4),
\newblock \emph{Journal of Combinatorial Theory, Series A} \textbf{25} (1978), 288--318.
\newblock DOI: \url{https://doi.org/10.1016/0097-3165(78)90021-3}.

\bibitem{MS1977}
F. J. MacWilliams and N. J. A. Sloane,
\newblock \emph{The Theory of Error-Correcting Codes},
\newblock North-Holland, Amsterdam, 1977.

\bibitem{MT2018}
J. V. S. Morales and H. Tanaka,
\newblock An Assmus--Mattson theorem for codes over commutative association schemes,
\newblock \emph{Designs, Codes and Cryptography} \textbf{86} (2018), 1039--1062.
\newblock DOI: \url{https://doi.org/10.1007/s10623-017-0376-y}.

\bibitem{NW2025}
H. Nozaki and Y. Watanabe,
\newblock Combinatorial characterzations of $T$-designs in the nonbinary Johnson scheme,
\newblock arXiv:2512.22034v1, 2025.
\newblock \url{https://arxiv.org/abs/2512.22034v1}.

\bibitem{Ozeki1997}
M. Ozeki,
\newblock On the notion of Jacobi polynomials for codes,
\newblock \emph{Mathematical Proceedings of the Cambridge Philosophical Society} \textbf{121} (1997), 15--30.
\newblock DOI: \url{https://doi.org/10.1017/S0305004196001399}.

\bibitem{TAG1985}
H. Tarnanen, M. J. Aaltonen, and J.-M. Goethals,
\newblock On the nonbinary Johnson scheme,
\newblock \emph{European Journal of Combinatorics} \textbf{6} (1985), 279--285.
\newblock DOI: \url{https://doi.org/10.1016/S0195-6698(85)80039-1}.

\bibitem{WZ2026}
X. Wu and J. Zhou,
\newblock Research on Linear Codes Holding $q$-Ary $t$-Designs,
\newblock arXiv:2603.12761v1, 2026.
\newblock \url{https://arxiv.org/abs/2603.12761v1}.

\end{thebibliography}
\end{document}